# Unsteady gas dynamics modeling for leakage detection in parallel pipelines

Ilgar G. Aliyev[*1], Konul A. Gafarbayli[2a],
Ahad J. Mammadov[2] and Mammadrzayeva Firangiz[2]

[1]*Head of Operation and Reconstruction of Buildings and Facilities Department,
Azerbaijan University Architecture and Construction, Baku, Azerbaijan*
[2]*Department of Land Reclamation and Water Resources Construction,
Azerbaijan University Architecture and Construction, Baku, Azerbaijan*

**Abstract.** This study presents a novel analytical framework for modeling unsteady gas dynamics in parallel pipeline systems under leakage conditions. The proposed method introduces a time-dependent leakage mass flow rate function, $G(t)=K \times P_1 \times e^{-\beta t}$, which dynamically captures the temporal decay of leakage based on real-time inlet pressure measurements. This functional form allows for a more physically consistent and mathematically tractable representation of gas loss compared to conventional constant-rate or stepwise models. The pipeline system is partitioned into three regions relative to the leakage point, and closed-form pressure solutions are derived using Laplace transform techniques. These expressions enable direct estimation of the leakage location through inverse pressure profiles, eliminating the need for computationally intensive iterative schemes. The analytical model is further validated against representative benchmark scenarios, demonstrating good agreement with literature-based results. A comparative analysis underscores the model's ability to localize leakage using minimal sensor data while preserving interpret ability-an essential feature for deployment in industrial environments. The approach provides a lightweight yet robust alternative to purely numerical or machine learning-based solutions and offers potential integration into real-time monitoring systems. This work contributes to the field by unifying gas dynamic principles, sensor-assisted modeling, and analytical solution strategies to enhance the reliability and speed of leak detection in modern gas transport infrastructures.

## 1. Introduction

Ensuring the integrity of gas pipelines is critical to operational safety and environmental protection. In parallel pipeline networks operating under unsteady flow conditions, the presence of leaks introduces significant challenges to system analysis and control. One of the key problems in this context is the incorporation of the time-varying leakage function $G(t)$ into the gas dynamic boundary conditions. Prior studies have commonly assumed $G(t)$ to be constant due to the indeterminacy introduced by its dependency on both the leak location and the actual mass flow rate of the escaping gas (Kowalsky *et al.* 2014, Souza *et al.* 2020).

However, gas leakage is highly sensitive to several physical factors, including the pressure gradient along the pipeline, the size and geometry of the leakage point, and whether the pipeline is installed above or below ground. It is known that gas leaks from high-pressure zones at higher flow rates compared to low-pressure zones, and larger leak openings result in higher mass flow. Therefore, any realistic model must incorporate spatially and temporally variable leakage behavior.

This paper proposes a new analytical expression for $G(t)$ based on pressure decay at the pipeline inlet. By replacing the unknown $G(t)$ with a measurable expression, $G(t)=f(P(0,t))$, and determining the coefficient via least squares fitting to experimental pressure data, we develop a model that accounts for leak location and severity without direct measurement of leak parameters. This approach is particularly feasible today thanks to the availability of IoT-based pressure sensors for real-time monitoring. Where *f* is a proportionality constant determined by system-specific parameters such as pipeline characteristics and gas properties.

Traditional detection methods rely on manual inspections or fixed threshold alarms, which are often delayed and inaccurate (Abhulimen *et al.* 2007, Brunone *et al.* 2015, Ellali *et al.* 2024). As pipeline networks become more complex and extensive, particularly with the use of parallel systems for redundancy and pressure balancing, the need for reliable, real-time detection models becomes essential.

This study focuses on the analytical modeling of unsteady gas dynamics in parallel pipelines under leakage scenarios. A key innovation of this work lies in handling the time-dependent leakage function, which has previously been assumed constant due to its mathematical complexity. We propose a dynamic and measurable alternative linked to observable pressure fluctuations, which can be effectively monitored using IoT technologies. Our model allows for the localization and quantification of leakages without requiring prior knowledge of leakage coordinates.

## 2. Related work

The problem of gas leakage detection and localization in pipeline systems has been addressed through various modeling strategies, which can be broadly classified into three principal categories:

### *2.1 Steady-state and pressure deviation approaches*

Conventional methods for leak detection often rely on pressure drop analysis under steady-state assumptions, where deviations from expected profiles are interpreted as indicative of leaks (Kowalsky *et al.* 2014, Souza *et al.* 2020, Brunone *et al.* 2015). While simple and computationally efficient, these methods lack the temporal resolution to capture transient leak behavior or detect small, evolving leaks.

### *2.2 Transient-based and inverse modeling techniques*

Table 1 Key distinguishing features of the study

| Aspect | Conventional Approaches | This Study |
|---|---|---|
| Leakage Function | Constant or predefined leak rate | Time-dependent exponential decay model G(t) |
| Solution Methodology | Numerical (iterative, discretized) | Closed-form analytical expressions using Laplace transform |
| Leak Localization | Inverse problem or supervised classification | Direct analytical estimation from pressure functions |
| Data Utilization | Extensive training data or synthetic benchmarks | Sensor-calibrated with minimal input |
| Applicability in Real Time | Often limited by computational cost | Suitable for online applications due to explicit expressions |

More advanced studies have considered unsteady flow dynamics by solving the governing equations of gas motion under leakage conditions, frequently employing inverse problem frameworks (Ellali *et al.* 2024, Samuel *et al.* 2024). Such models are generally solved numerically using finite difference or finite volume schemes and are capable of more precise leak localization. However, their reliance on iterative solution procedures can limit real-time applicability.

### 2.3 Data-driven and AI-assisted methods

In recent years, machine learning-based methods have gained prominence in pipeline monitoring. Techniques such as Long Short-Term Memory (LSTM) networks, Convolutional Neural Networks (CNN), and hybrid AI models have been successfully applied to Supervisory Control and Data Acquisition (SCADA) data for real-time leakage classification and location prediction (Yongmei *et al.* 2021, Firouzi *et al.*2021, Aliyev 2024, Turkach *et al.* 2023). Despite their high accuracy, these approaches typically require extensive historical datasets for training and may suffer from interpretability limitations in safety-critical contexts.

### 2.4 Analytical framework of this study

The proposed analytical model introduces a novel leakage function $G(t)=K \times P_{in}(t) \times e^{-\beta t}$, representing a time-decaying mass flow rate of leaked gas, dynamically linked to inlet pressure. This formulation contrasts sharply with prior works that assume a constant or stepwise leakage rate. The main differentiating features of the present study are outlined in Table 1.

The proposed framework thus contributes a mathematically robust and practically efficient method for leakage detection in parallel gas pipelines. The method leverages both theoretical insights and sensor-based adaptability, positioning it as a viable alternative to purely empirical or computationally intensive approaches.

The study aligns with the growing trend of combining physics-based models with sensor-informed calibration for real-time infrastructure diagnostics. For instance, (Li *et al.* 2024) developed a hybrid LSTM-Transformer network achieving 99.995% accuracy in leak prediction over long distances, while (Xuguang *et al.* 2024) addressed imbalance in data samples through multi-scale kernel temporal convolutional networks (MKTCN). Our approach complements these methods by offering analytical transparency and computational simplicity, which are advantageous in regulatory and industrial environments where interpretability is crucial.

Moreover, recent reviews emphasize the potential of digital twin models for leakage forecasting, yet highlight the need for mathematically tractable core models for integration (Zhang *et al.* 2024, Aliyev 2024, Kamel *et al.* 2023). The current study provides such a model and demonstrates its practical use through closed-form analysis validated against representative scenarios.

The pressure decay behavior at the inlet of a gas pipeline is not uniform; it varies significantly depending on the location, size, and shape of the leak, as well as the mass flow rate of the escaping gas. Specifically, when the leak occurs closer to the inlet, the rate of pressure drop at the inlet tends to be higher. Conversely, when the leak is located near the end of the pipeline, the pressure decline at the inlet is less pronounced.

Despite these variations, extensive experimental studies have shown that the overall pressure decay at the inlet follows an exponential trend over time. For instance, (Shi *et al.* 2021) demonstrated that leakage-induced pressure signals attenuate exponentially with time under a range of operational and physical scenarios.

Therefore, the expression $P(0,t)=Pe^{-\beta t}$ used in this study to model inlet pressure dynamics is not only mathematically convenient but also physically realistic. It provides a reliable basis for defining the time-dependent leakage mass flow rate function $G(t)$, which is linked to inlet pressure via.

## 3. Methodology

The methodology of this study integrates analytical modeling, Laplace transform techniques, and empirical parameter estimation to investigate unsteady gas dynamics in parallel pipeline systems under leakage conditions. The overall approach consists of three interconnected phases: mathematical modeling, functional approximation using experimental data, and analytical solution derivation.

### 3.1 System configuration and modeling assumptions

We consider a parallel gas pipeline system operating under a unified hydraulic regime, where both pipelines share common inlet and outlet boundaries. In the event of a leakage, the system is divided into three segments. Each segment is governed by a one-dimensional unsteady gas dynamics equation.

Key assumptions include:
- The flow is compressible and unsteady, but one-dimensional.
- Initial and boundary pressures are known.
- Leak is modeled as a point-source discontinuity at location $\ell$.
- The pressure field is symmetric and decays exponentially over time due to leakage.

### 3.2 Governing equations

To address the nonlinear nature of the momentum equation—particularly the term involving $v|v|$—a Chernov-type linearization is applied to facilitate the solution:
$$v|v| \approx |v_0|v$$
where $v_0$ is a known or previously estimated velocity profile. This approximation transforms the nonlinear term into a quasi-linear form, enabling the use of stable and efficient numerical schemes. The resulting linearized form of the friction term becomes:

$$2a = \lambda v_0 / 2d$$

This linearization simplifies the governing equations while preserving the essential dynamics of frictional losses. It allows for tractable implementation in iterative solvers or semi-implicit schemes, making the model suitable for real-time simulation and control applications.

### *3.3 Application of laplace transform*

The Laplace transform is applied to the system of governing partial differential equations to facilitate the handling of time-dependent boundary conditions and to linearize the unsteady dynamics. Transformed equations allow for the incorporation of localized leakage effects using delta-function modeling at $x=\ell$. Boundary and initial conditions are also transformed accordingly.

The solution in the Laplace domain is then inverted using inverse transform techniques to yield closed-form expressions for the pressure functions $P_1(x,t)$, $P_2(x,t)$, and $P_3(x,t)$.

### *3.4 Leakage function approximation and calibration*

The leakage mass flow rate $G(t)$, which is difficult to determine analytically due to dependence on the unknown leak location, is approximated using inlet pressure data. The decay coefficient $\beta$ is estimated using the least squares method from real or simulated pressure data collected via IoT-based wireless pressure sensors.

### *3.5 Scenario-based validation*

To validate the model, leakage scenarios are simulated with leak locations (25 km). Using predefined pipeline characteristics and experimentally derived values for $\beta$, the pressure distribution is computed and analyzed along the pipeline axis. The results are compared in terms of spatial response and dynamic behavior.

### *3.6 Numerical implementation*

For practical computation:
- A truncated Fourier series (first 20 terms) is used to approximate the infinite series in the pressure solution.
- The model is evaluated at selected time intervals (e.g., $t=120$ *c*) for comparison across scenarios.
- Graphical visualization is used to interpret the pressure field evolution and validate the model's sensitivity to leak position.

## 4. Theoretical model

We consider a parallel gas pipeline system operating under identical inlet and outlet pressures, thereby forming a closed-loop network. In the presence of a leak in one of the pipelines, the system is conceptually divided into three distinct regions for the purposes of unsteady gas

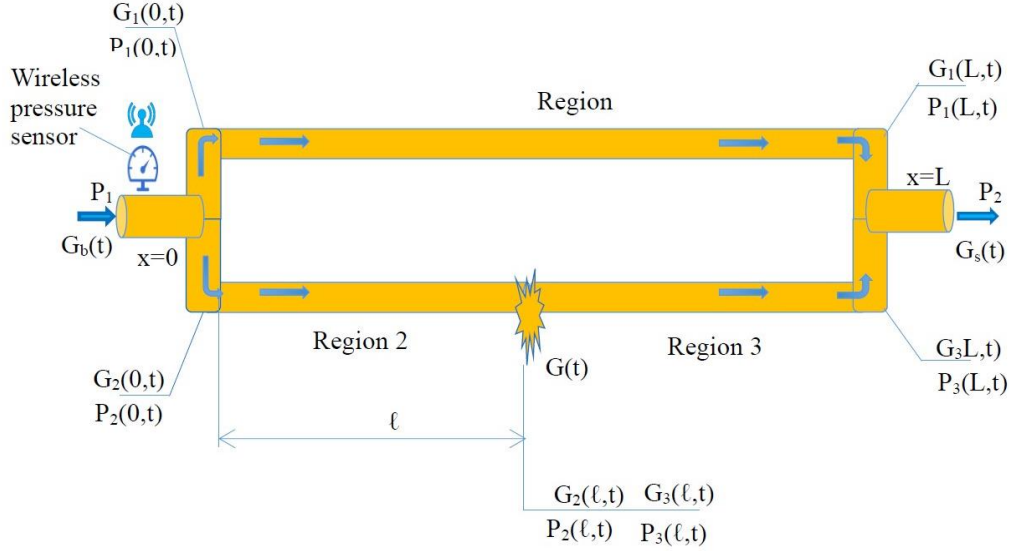

Fig. 1 Parallel gas pipeline system operating under identical inlet and outlet pressures

dynamics modeling (Fig. 1).

Region 1: The intact pipeline segment through which gas flows without leakage.
Region 2: The segment of the defective pipeline extending from the inlet up to the leakage point.
Region 3: The segment of the defective pipeline from the leakage point to the outlet.

Each of these regions is governed by the one-dimensional unsteady gas dynamics equations, yielding a set of pressure functions $P_1(x,t)$, $P_2(x,t)$, and $P_3(x,t)$, and corresponding mass flow rates $G_1(t)$, $G_2(t)$, and $G_3(t)$. This results in six unknown functions, necessitating six independent boundary and interface conditions to close the system.

### 4.1 Boundary and matching conditions

Four of the required conditions are derived from the continuity of pressure and conservation of mass flow at the pipeline inlet and outlet:

At the Inlet ($x=0$):
$$P_1(0,t)=P_2(0,t), \quad G_1(t)+G_2(t)=G_0(t)$$

At the Outlet ($x=L$):
$$P_1(L,t)=P_3(L,t), \quad G_1(t)+G_3(t)=G_k(t)$$

The remaining two conditions are applied at the leakage location $x=\ell$, assuming the leak behaves as a point discontinuity:

At the Leak Point ($x=\ell$):
$$P_2(\ell,t)=P_3(\ell,t), \quad G_2(t)=G_3(t)+G(t)$$

Here, $G(t)$ denotes the time-dependent mass flow rate of the leaking gas, and the Heaviside function $H(x-\ell)$ is used to localize the leakage at point $\ell$. This function equals 0 for $x<\ell$ and 1 for $x\geq\ell$, ensuring that leakage is modeled only downstream of the rupture location.

*4.2 Modeling challenges and leakage function approximation*

A central difficulty in the modeling process lies in accurately characterizing the leakage term $G(t)$, which is both time-dependent and spatially localized. Since the leak location $\ell$ is typically unknown, the function $G(t)$ cannot be determined analytically. Moreover, the leakage rate depends not only on the pressure but also on the rupture geometry, pipe burial depth, and local flow regime. This makes it impractical to define $G(t)$ without incorporating empirical or sensor-based data.

To overcome this, we approximate the leakage mass flow rate using the measurable inlet pressure $P_1(0,t)$. Empirical observations show that inlet pressure decays faster when leakage occurs closer to the pipeline inlet or when the leak is larger. This motivates the approximation:

The key innovation is replacing $G(t)$ by $P_1(0,t)$, reducing the number of unknowns and linking the problem to directly measurable inlet pressure data. Assuming exponential pressure decay: Based on empirical observations and experimental studies, we propose the following approximation:

$$G(t) \approx (g P_1(0,t))/(2ac^2) \, H(x-\ell)$$

where:
- $g$ is the gravitational acceleration, m/s$^2$
- $2a$ Chernov-type linearization coefficient, s$^{-1}$
- $c$ is the speed of sound in the gas medium, m/s

This form allows us to express $G(t)$ as a function of directly measurable quantities, thereby reducing the number of unknowns and enabling real-time leak modeling.

Assuming that the inlet pressure decays exponentially during leakage events (Shi *et al.* 2021)

$$P_1(0,t) = P_1 \, e^{-\beta t} \qquad (1)$$

we obtain:

$$G(t) = (g \, P_1)/(2ac^2) \, e^{-\beta t} \, H(x-\ell)$$

*4.3 Calibration via initial conditions*

However, this expression may not match the actual initial leakage rate $G_0 = G(0)$. To align the model with known system behavior at $t=0$, we introduce a correction factor $K$, leading to the modified leakage function

$$G(t) = (K \, g \, P_1)/(2ac^2) \, e^{-\beta t} H(x-\ell) \qquad (2)$$

This expression serves as a physically motivated approximation for transient leakage flow originating at location $x=\ell$. Below is a breakdown of the terms:

$K$- A leakage coefficient that incorporates geometric and material properties of the leakage point.

$g$- Gravitational acceleration, included to account for hydrostatic pressure contributions.

$P_1$- The inlet pressure, providing the initial driving force for gas flow through the leak.

$a$-Represents the Chernov linearization factor.

$\beta$-Exponential decay constant, governing how quickly the leakage rate declines due to pressure reduction over time.

$H(x-\ell)$- Heaviside step function, activating the leakage term only at the spatial location $x=\ell$.

The form of the function reflects realistic transient leakage behavior: the leakage initiates at

*x=ℓ* and decays over time as the internal pressure drops. The exponential term *e−βt* approximates the dissipation of leakage mass flow due to decreasing pressure and stabilizing flow conditions. The use of Charni's linearization for a ensures that internal friction and pipe geometry are properly accounted for, making the leakage flow rate responsive to pipeline design.

This leakage function acts as a source term in the pressure wave equation and directly affects the solution *P(x,t)*, as shown in the derivations in Section 4.

The coefficient *K* is calibrated by enforcing:

Using Eq. (1), the coefficient *K* at time *t*=0 is determined as follows.

$$G(0)=G_0 \Rightarrow K=(2ac^2 G_0)/(g P_1)$$

This ensures that the leakage model is consistent with the observed or expected initial conditions. Once system parameters such as $G_0$, $P_1$, *a*, and *c* are known, *K* can be computed at commissioning and used for dynamic simulation and fault diagnostics.

### *4.4 IoT-Enhanced pressure monitoring and leak localization*

Modern IoT technologies enable the continuous collection of pressure data via wireless sensors installed along the pipeline (Ekong *et al.* 2024, Effiom *et al.* 2024). These sensors transmit real-time measurements to a control center, allowing early detection of abnormal pressure drops and facilitating leak localization.

The *β* coefficient is determined via least squares fitting from pressure data collected by IoT sensors. Assuming exponential decay: The βcoefficient can be extracted by linear regression on the logarithmic form of pressure data. IoT-Based Pressure Monitoring Modern IoT technology enables real-time acquisition of pressure data using wireless sensors. These sensors transmit measurements to a control center over the Internet, allowing remote monitoring of pipeline conditions. The integration of IoT significantly enhances the accuracy of leakage detection and parameter estimation.

## 5. Mathematical model

The unsteady gas flow in the parallel pipeline system is modeled using the linearized one-dimensional gas dynamics equations. For each segment *i*=1,2,3 of the network, the system of governing equations includes the continuity and momentum conservation laws. The equations are expressed in terms of pressure $P_i(x,t)$ and mass flux $G_i(x,t)$ as follows: The governing equations are derived from the principles of mass and momentum conservation.

In the system layout presented in the Fig. 1, a pipeline with a rupture is divided into three segments:

Segment 1-the intact branch of the network ($0 \leq x \leq L$);
Segment 2-the section from the junction to the leak point ($0 \leq x \leq \ell$);
Segment 3-the section from the leak point to the end ($\ell \leq x \leq L$).
Governing equations (Aliyev 2024)

$$-\frac{1}{c^2}\frac{\partial P_i(x,t)}{\partial t}=\frac{\partial G_i(x,t)}{\partial x} \tag{3}$$

$$-\frac{\partial P_i(x,t)}{\partial x} = 2aG_i(x,t) \tag{4}$$

Taking the derivative of Eq. (3) with respect to x, and substituting the resulting expression for $\frac{\partial G_i(x,t)}{\partial x}$ into Eq. (2), we derive the following diffusion-type equation for the pressure $P_i(x,t)$

$$\frac{\partial^2 P_i(x,t)}{\partial x^2} = \frac{2a}{c^2}\frac{\partial P_i(x,t)}{\partial t}, \quad i=1.2.3 \tag{5}$$

### 5.1 Boundary and initial conditions

To close the problem, we define the boundary and initial conditions reflecting the pipeline configuration:

At inlet ($x=0$):

$$\begin{cases} P_1(x,t) = P_2(x,t) \\ \dfrac{\partial P_1(x,t)}{\partial X} + \dfrac{\partial P_2(x,t)}{\partial X} = -2aG_0(t) \end{cases}$$

At outlet ($x=L$):

$$\begin{cases} P_1(x,t) = P_3(x,t) \\ \dfrac{\partial P_1(x,t)}{\partial X} + \dfrac{\partial P_2(x,t)}{\partial X} = -2aG_k(t) \end{cases}$$

Leak interface ($x=\ell$):

$$\begin{cases} P_2(x,t) = P_3(x,t) \\ \dfrac{\partial P_3(x,t)}{\partial X} - \dfrac{\partial P_2(x,t)}{\partial X} = -2aG(t) = -K\dfrac{gP_1}{c^2}e^{-\beta t} \end{cases}$$

Initial pressure distribution:
$$t=0;\ P_i(x,0) = P_1 - 2aG_ix:\ i=1,2,3$$

### 5.2 Boundary behavior and early-time approximation

In the study of non-stationary gas dynamics in parallel pipelines operating under a unified hydraulic regime, the formulation of appropriate boundary conditions is essential for both analytical treatment and numerical implementation of the model. When leakage occurs in one of the parallel pipelines, the total inlet and outlet mass flow rates, denoted $G_0(t)$ and $G_k(t)$, become time-dependent and may deviate from each other due to the local mass imbalance introduced by the leakage.

However, an important simplification can be employed during the early stage following the onset of leakage. Because the pressure disturbance caused by the leakage propagates at the finite acoustic velocity $c$ of the gas (typically $c \approx 383$ m/sec), its effect at the pipeline boundaries is not instantaneous. For a leakage located at a distance $\ell$ from the boundaries, the delay time $\tau$ for the pressure wave to reach the inlet and outlet can be estimated as:

$$\tau = \frac{\ell}{c}$$

For instance, if the leakage is located at $\ell=50$ km, then:
$$\tau = \frac{50000}{383} \approx 130 \text{ seconds}$$

Within this transient time interval $t\in[0,\tau)$, the difference between inlet and outlet flow rates is physically negligible. Therefore, it is justified to assume:
$$G_b(t)=G_s(t), \quad \text{for } 0 \leq t<\tau.$$

This assumption allows for a significant simplification of the governing system of equations, avoiding complex coupling through time-varying boundary conditions and facilitating a more stable and computationally efficient solution approach.This modeling strategy offers a balance between physical fidelity and computational tractability, making it suitable for practical applications such as real-time leakage detection and control in large-scale pipeline systems.

The application of the Laplace transform allows for a more stable and analytically tractable representation of pressure functions in unsteady gas dynamic systems. This method enables a more accurate characterization of the system's dynamic response by incorporating both boundary and initial conditions, particularly facilitating the modeling of local disturbances such as leakage.

### 5.3 Mathematical model innovation

A key innovation of the present study lies in the introduction of a time-dependent leakage mass flow rate function $G(t)$, into the unsteady gas dynamics model of parallel pipelines. Unlike previous studies that treat leakage as a constant loss or impose fixed boundary anomalies, this work formulates $G(t)$ as a dynamically evolving function governed by real-time pressure decay, which better reflects the physical behavior of gas leaks in operational systems.

The function is modeled as an exponentially decaying expression:

$G(t)=K\times P_1\times e^{-\beta t}$, where the coefficients $K$ and $\beta$ are determined using pipeline specifications and least-squares fitting of IoT-based sensor data, respectively. This approach enables direct analytical linkage between pressure dynamics and the leak's location and magnitude.

The originality of this method is twofold:

It provides a closed-form solution for the pressure distribution that includes a non-constant source term ($G(t)$), thus enabling analytical inversion to detect leakage coordinates. It introduces a hybrid analytical–data-driven paradigm where sensor-informed calibration enhances the reliability of theoretical predictions.

To the best of the authors' knowledge, no existing analytical models incorporate a time-dependent leakage mass flow rate in such a form for inverse problem solving in gas pipeline systems. This contributes significantly to both the theoretical foundation of unsteady flow modeling and practical leakage detection methods.

Comparative Context with Prior Works

Previous works in the field of gas leakage modeling, such as those by (Petro *et al.* 2020, Wang *et al.* 2018, Fang *et al.* 2022), typically employ one of the following approaches:

-Constant Leakage Assumption: Leakage is modeled as a constant mass outflow, independent of system pressure or time.

-Static Pressure Drop Localization: Leak detection is based on steady-state deviation of pressure measurements at fixed nodes.

| Approach Type | Leakage Function $G(t)$ | Analytical Solution | Real-time Calibration | Inverse Problem Capability |
|---|---|---|---|---|
| Constant Leakage Models | $G(t)$=const | Limited | No | No |
| Static Deviation Localization | Not explicitly modeled | No | No | Partial |
| Purely Numerical Inverse | Variable, simulation-based | No | Possible | Yes (but no closed-form) |
| **This Study** | $G(t)=K \times P_{in}(t) \times e^{-\beta t}$ | **Yes** | **Yes (IoT-based)** | **Yes (Analytical)** |

-Numerical-only Strategies: Fully numerical inverse models that require complex iterative solvers and lack analytical transparency.

As shown above, the present approach uniquely offers a closed-form analytical solution to the inverse problem while incorporating real-time variability via the calibrated function $G(t)$. This enhances interpretability, reduces computational burden, and bridges the gap between theory and practice in leakage detection.

## 5.4 Physical ınterpretation and practical ımplications

- At early times ($t \to 0$), the exponential terms can be approximated using Taylor expansions, revealing a nearly uniform response dominated by the initial gas release.
- At long times ($t \to \infty$), both exponential terms decay to zero, and the system tends toward a quasi-stationary state governed by $1/\beta$, indicating the asymptotic behavior of the leakage impact.
- The Fourier structure makes it possible to analyze localized disturbances caused by leakage, particularly useful for leak localization.

This pressure function is particularly well-suited for real-time diagnostics, especially when coupled with IoT-based pressure monitoring systems. The presence of measurable parameters such as $P_1$, $\beta$, and $K$, combined with rapid convergence of the series expansion, makes the model both physically realistic and computationally efficient.

It should be noted that in the derived pressure expression, the only unknown parameter is the exponential decay coefficient $\beta$.

As discussed in the Theoretical Model section of the paper,

this coefficient can be estimated using pressure data collected in real-time via wireless sensors enabled by Internet of Things (IoT) technology. Applying the least squares method to the sensor data allows accurate determination of $\beta$.

Naturally, the sooner this estimation is performed after a leakage event, the more accurate the value of $\beta$ will be.

Additionally, as mentioned in the paper, the approximation $G_b(t) \approx G_s(t)$ holds valid for small time intervals immediately following the onset of leakage. All other parameters in the pressure equation are known from the technical specifications (passport) of the given gas pipeline, which allows precise calculation of the calibration coefficient $K$.

## 6. Theoretical–experimental analysis

To evaluate the applicability and accuracy of the derived pressure distribution model under real-world conditions, a practical leakage scenario is considered in a parallel gas pipeline system operating under a unified hydraulic regime (*A*8, *A*9, *A*10). The objective of the analysis is to assess the variation of gas pressure along the pipeline axis as a function of the leakage location. In this case, the leakage is assumed to occur at a distance of $\ell=25$ km from the pipeline inlet.

Table 2 Time-dependent variation of inlet pressure $P(0,t)$ due to loss of integrity at $\ell=25$ km in a parallel gas pipeline (Aliyev 2024)

| Time $t$ (s) | Pressure $P(0,t)$ [$\times 10^4$ Pa] |
|---|---|
| 0 | 55.0 |
| 300 | 42.38 |
| 600 | 40.37 |
| 900 | 38.3 |
| 1200 | 36.21 |
| 1500 | 34.11 |

The parameter values were chosen based on typical operational ranges in medium-pressure natural gas transmission pipelines. Specifically,

$P_1=5.5\times10^5$ Pa; $G_0=30$ Pa·s/m; $2a=0.1$ s$^{-1}$: $c=383.3$ m/s; $L=100$ km; $g=9.81$ m/s²

Step 1: Determination of Calibration Coefficient $K$

Using the known parameters and the initial mass flow rate $G_0$, the calibration coefficient $K$ in the pressure equation is computed using the formula:

$$K=(2a\ c^2\ G_0) / (g\ P_1)=0.802$$

This coefficient ensures that the leakage function $G(t)$ starts from the correct initial value and is consistent with the physical configuration of the pipeline.

Step 2: Estimation of Pressure Decay Coefficient $\beta$ Using Experimental Data

To model the pressure decay function:

$$P(0,t)=P_1 e^{-\beta t}$$

we utilize pressure-time data obtained from previous studies that used pipelines with identical physical characteristics and constant leakage rates. Separate datasets are used for each scenario, corresponding to leaks at the respective positions.

The coefficient $\beta$ is estimated by applying the least squares fitting method to the logarithmic form of the measured pressure data:

According to reference (Aliyev 2024), the pressure-time characteristics at the pipeline inlet were evaluated under the condition that the leakage mass flow rate is known and fixed as $G(t)=0.8G_0$. Specifically, in a gas pipeline network operating under a unified hydraulic regime, when the pipeline integrity is compromised at a distance of $\ell=25$ km, the inlet pressure decreases from $P(0,0)=55\times10^4$ Pa to $P(0,300)=42.38\times10^4$ Pa over a 300-second interval.

Based on these pressure measurements, and applying the least squares fitting method to the exponential pressure decay model, the decay coefficients $\beta$ corresponding to leakage scenario are determined. To determine a physically realistic value for $\beta$, we utilize theoretical pressure decay data presented in reference (Aliyev 2024), where the following pressure values were observed at the pipeline inlet x=0 under a leakage condition located at $\ell=25$ km:

Assuming that the pressure decay follows an exponential law of the form:

Table 3 Pressure distribution along the pipeline for each segment

| x, km | t=100 sec | | | t=300 sec | | | t=600 sec | | |
|---|---|---|---|---|---|---|---|---|---|
| | $P_1(x,t)$ [×10⁴ Pa] | $P_2(x,t)$ [×10⁴ Pa] | $P_3(x,t)$ [×10⁴ Pa] | $P_1(x,t)$ [×10⁴ Pa] | $P_2(x,t)$ [×10⁴ Pa] | $P_3(x,t)$ [×10⁴ Pa] | $P_1(x,t)$ [×10⁴ Pa] | $P_2(x,t)$ [×10⁴ Pa] | $P_3(x,t)$ [×10⁴ Pa] |
| 0 | 53.40 | 53.40 | - | 45.39 | 45.39 | - | 34.60 | 34.60 | 0.00 |
| 5 | 53.41 | 51.21 | - | 47.16 | 41.35 | - | 37.44 | 29.63 | 0.00 |
| 10 | 53.15 | 48.22 | - | 48.46 | 36.61 | - | 39.87 | 24.14 | 0.00 |
| 15 | 52.55 | 44.07 | - | 49.15 | 31.02 | - | 41.71 | 18.03 | 0.00 |
| 20 | 51.96 | 38.31 | - | 49.58 | 24.37 | - | 43.25 | 11.14 | 0.00 |
| 25 | 51.19 | 32.69 | 32.69 | 49.57 | 18.49 | 18.49 | 44.29 | 5.28 | 5.28 |
| 30 | 50.52 | - | 36.81 | 49.45 | - | 22.87 | 45.11 | 0.00 | 9.64 |
| 35 | 49.71 | - | 41.07 | 49.02 | - | 28.02 | 45.53 | 0.00 | 15.03 |
| 40 | 49.03 | - | 43.72 | 48.60 | - | 32.11 | 45.82 | 0.00 | 19.64 |
| 45 | 48.22 | - | 45.21 | 47.95 | - | 35.35 | 45.77 | 0.00 | 23.63 |
| 50 | 47.53 | - | 45.90 | 47.37 | - | 37.89 | 45.68 | 0.00 | 27.10 |
| 55 | 46.72 | - | 45.91 | 46.62 | - | 39.66 | 45.31 | 0.00 | 29.94 |
| 60 | 46.03 | - | 45.65 | 45.98 | - | 40.96 | 44.95 | 0.00 | 32.37 |
| 65 | 45.21 | - | 45.05 | 45.18 | - | 41.65 | 44.35 | 0.00 | 34.21 |
| 70 | 44.54 | - | 44.46 | 44.51 | - | 42.08 | 43.79 | 0.00 | 35.75 |
| 75 | 43.71 | - | 43.69 | 43.70 | - | 42.07 | 43.01 | 0.00 | 36.79 |
| 80 | 43.04 | - | 43.02 | 43.01 | - | 41.95 | 42.29 | 0.00 | 37.61 |
| 85 | 42.21 | - | 42.21 | 42.18 | - | 41.52 | 41.35 | 0.00 | 38.03 |
| 90 | 41.53 | - | 41.53 | 41.48 | - | 41.10 | 40.45 | 0.00 | 38.32 |
| 95 | 40.72 | - | 40.72 | 40.62 | - | 40.45 | 39.31 | 0.00 | 38.27 |
| 100 | 40.03 | - | 40.03 | 39.87 | - | 39.87 | 38.18 | 0.00 | 38.18 |

$$P(0,t) = P_1 \, e^{-\beta t}$$

we apply a least squares fitting procedure using these data points. The optimal value of $\beta$ is obtained by minimizing the error between the observed pressures and the modeled exponential curve. Using Python's curve fitting method, we find that the best-fitting value is:

$$\beta = 1.03 \times 10^{-4} \text{ s}^{-1}$$

This value of $\beta$ closely reproduces the observed pressure at $t=300$ s, where $P(0,300) \approx 42.38 \times 10^4$ Pa, as recorded in the baseline reference.

With the values of $K$, $\beta$, and other parameters determined, we employ the analytical expression for pressure distribution:

Using the estimated decay coefficient $\beta$ for the leakage location at $\ell=25$ km, analytical pressure distributions $P_1(x,t)$, $P_2(x,t)$ and $P_3(x,t)$ along the parallel gas pipeline are computed for all three segments. The calculations are performed for the time range from $t=0$ to $t=600$ seconds with increments of 100 seconds, and for the spatial domain from $x=0$ to $x=100$ km with increments of 5 km.

The results of the pressure distributions across the pipeline for each segment are presented in the table below ($t=100$, $t=300$, $t=600$).

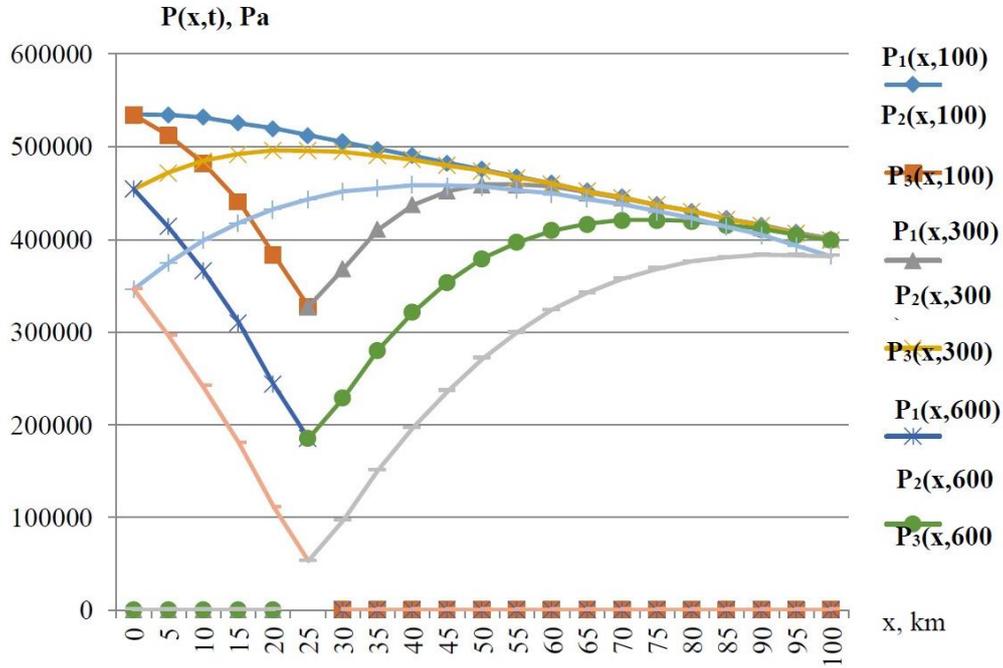

Fig. 2 Time-dependent pressure distribution in three zones of the damaged gas pipeline ($\ell$=25 km)

*6.1 Model accuracy and comparative analysis*

At $t$=300 seconds, the computed pressure at the pipeline inlet using the proposed model with a leakage at $\ell$=25 km yields $P(0,300)=45.39\times10^4$ Pa. In comparison, the theoretical reference (Aliyev 2024), which assumes a constant leakage rate $G(t)=G_0$, provides $P(0,300)=42.38\times10^4$ Pa. The relative deviation between these two results is approximately 7%, which falls within an acceptable range for engineering modeling.

Two main factors contribute to this discrepancy:

1. The value of the decay coefficient $\beta$ is estimated from available pressure data and may contain a degree of fitting uncertainty.

2. More fundamentally, the leakage function $G(t)$ in the proposed model is pressure-dependent and exhibits exponential decay over time. This dynamic relationship makes the system more responsive to pressure changes, unlike traditional models where $G(t)$ is held constant (Liang 2020, Suvalija *et al.* 2022, Turkach *et al.* 2023).

In the proposed model, as time progresses, the decreasing pressure leads to a corresponding reduction in the leakage rate $G(t)$. This dynamic interaction causes a smoother and less abrupt pressure drop at the inlet, aligning more closely with real-world pipeline behavior under leakage scenarios. The model thus captures the essential physics of pressure-dependent leakage and provides a more flexible and accurate predictive framework.

A joint analysis of Fig. 2 and Table 3 reveals several important insights. First, during the initial 200 seconds, the pressure at the end of the pipeline shows almost no change due to the leakage. Specifically, the pressure at $x=L$ remains approximately $P(L,300)=39.9\times10^4$ Pa, which is nearly equal to the initial pressure $P(L,0)=40\times10^4$ Pa.

This behavior is consistent with the explanation given in the Boundary Behavior and Early-Time Response section: the pressure wave propagates through the pipeline at the speed of sound in gas. Given that the leakage occurs at $\ell$=25 km, the remaining distance to the outlet is 75 km. With a sound speed of approximately 383.3 m/s, the pressure wave requires:

$$\tau = 75000/383.3 \approx 196 \text{ seconds}$$

to reach the outlet. Accordingly, for the first ~200 seconds, the outlet region is unaffected by the leak. It is only after this delay that a small pressure drop of $\Delta P \approx 0.1 \times 10^4$ Pa is observed over the next 100 seconds. By $t$=600 s, this drop increases to approximately $1.8 \times 10^4$ Pa.

In contrast, the inlet region shows a significant drop in pressure over the same time frame. At $x$=0, the pressure decreases from $P(0,0)=55 \times 10^4$ Pa to $P(0,600)=34.6 \times 10^4$ Pa, yielding a pressure drop of $20.4 \times 10^4$ Pa, which is approximately 11 times greater than the corresponding value at the outlet.

These observations confirm that the location of the leakage has a substantial impact on the pressure dynamics at the pipeline boundaries.

The result supports the earlier analytical conclusions regarding the non-uniform, unsteady gas flow behavior in a parallel pipeline system under leakage conditions and validates the physical foundations of the proposed dynamic gas flow model.

A key observation from Fig. 2 is that the maximum pressure drop consistently occurs at the leakage point. This confirms the sensitivity of the pressure distribution to the position of the leak and reinforces the core principle of the proposed model: the dynamic response of the gas system is predominantly governed by the leakage location within the parallel pipeline structure.

The analytical solution accurately reflects the transient nature of the pressure field, capturing both global decay and local spatial disturbances due to leakage. These results validate the model's effectiveness and its ability to adapt to different leakage locations without reparameterization.

Before presenting the leakage localization procedure, it is important to note that the underlying mathematical model and the analytical formulation of the pressure functions for unsteady gas dynamics under leakage conditions have been successfully validated. Based on these well-established expressions, the next goal is to determine the leakage location $\ell$ using an analytical-inverse approach applied to the pressure function at the pipeline inlet.

### 6.2 Leakage localization

The inlet pressure function $P_2(0,t)$, which depends on the leak location $\ell$, is matched with the sensor-based exponential pressure decay $P_1(0,t)=P_1 e^{-\beta t}$. This assumption allows the leakage point to be estimated using inverse modeling of pressure behavior at the pipeline inlet.

We expand the pressure function $P_2(0,t)$ using a Fourier series representation, isolating the dependence on $\ell$ within the cosine terms

$$P_2(0,t) = P_1 - \frac{KgP_1}{4aL} \times$$
$$\times \left[ \frac{1-e^{-\beta t}}{\beta} + 2\sum_{n=1}^{\infty}(-1)^n \cos\frac{\pi n(L-\ell)}{L} \frac{e^{-\beta t}-e^{-\alpha n^2 t}}{\alpha n^2 - \beta} \right] \tag{6}$$

Using the identity $\cos\frac{\pi n(L-\ell)}{L} = (-1)^n \cos\frac{\pi n \ell}{L}$, and dividing both sides of the equation by $P_1$, we obtain:

$$e^{-\beta t} = 1 - \frac{Kg}{4aL}\frac{1-e^{-\beta t}}{\beta} - \frac{Kg}{2aL}\sum_{n=1}^{\infty}\cos\frac{\pi n \ell}{L}\frac{e^{-\beta t}-e^{-\alpha n^2 t}}{\alpha n^2 - \beta} \tag{7}$$

By isolating the summation and expressing the right-hand side independently, we derive the final equation used in the numerical inversion

$$\sum_{n=1}^{\infty}\cos\frac{\pi n \ell}{L}\frac{e^{-\beta t}-e^{-\alpha n^2 t}}{\alpha n^2 - \beta} \approx \frac{2aL}{Kg}\left[1 - \frac{1-e^{-\beta t}}{4\beta} - e^{-\beta t}\right] \tag{8}$$

*6.3 Conclusion and validation of leakage localization method*

The results obtained through numerical analysis confirm the validity of the analytical formulation proposed for leakage localization in a parallel pipeline system. The equation used for estimating the leakage point $\ell$ integrates pressure function behavior and time-dependent decay, derived from real-time inlet sensor measurements.

To evaluate the accuracy of the proposed method, we used the derived inverse Eq. (8):
With the following input parameters:

$P_1 = 5.5 \times 10^5$ Pa
$G_0 = 30$ Pa·s/m
$a = 0.05$ s$^{-1}$
$c = 383.3$ m/s
$L = 100$ km
$g = 9.81$ m/s²
$K = 0.802$
$\beta = 1.03 \times 10^{-4}$ s$^{-1}$

In the leakage localization model, accurate determination of the leak location $\ell$ is theoretically validated by the convergence of the (17) left-hand side (LHS) and right-hand side (RHS) of the governing equation:

However, a direct comparison of these expressions reveals that the LHS and RHS possess different magnitudes and dimensional units. To evaluate the functional similarity and trend alignment between both sides, normalization is applied (Conesa *et al.* 2023). This technique rescales both LHS and RHS by their respective maximum absolute values:

Normalized LHS(*t*)=LHS(*t*)/max(|LHS(*t*)|), Normalized RHS(*t*)=RHS(*t*)/max(|RHS(*t*)|)

Through normalization, the analysis focuses on the temporal evolution and behavioral similarity of both expressions, rather than their absolute magnitudes. Although the original values at $t \approx 120$ s may differ in absolute magnitude (e.g., LHS≈5.99 and RHS≈37,880), the normalized plot shows that both expressions exhibit parallel trends and intersect in the same time range. This indicates a strong structural consistency between the two formulations, thereby reinforcing the model's theoretical validity and its practical capability to estimate the leakage location $\ell$.

In practical gas pipeline systems, the physical processes associated with leakage and pressure propagation are inherently unsteady and time-dependent. Eq. (8), which equates the summation-based LHS and analytical RHS, represents the core of the leakage localization model:

LHS=RHS

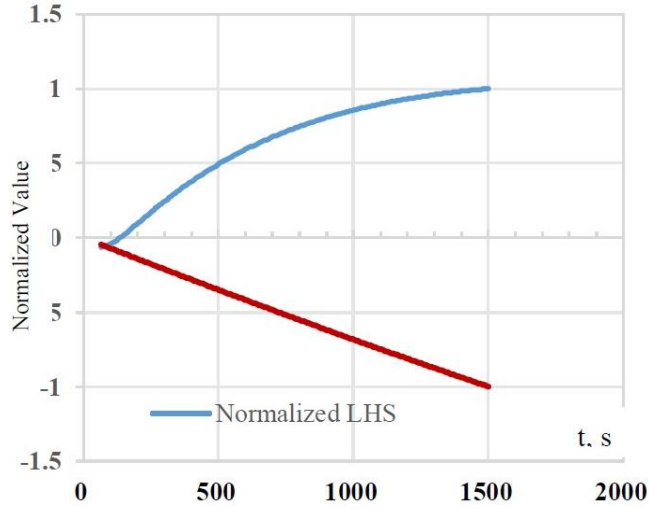

Fig. 3 Comparison of LHS and RHS in leakage localization equation

where the left-hand side is defined as:

$$\sum_{n=1}^{\infty} \cos\frac{\pi n \ell}{L} \frac{e^{-\beta t} - e^{-\alpha n^2 t}}{\alpha n^2 - \beta} \quad \text{and the right-hand side as:}$$

$$\frac{2aL}{Kg}\left[1 - \frac{1-e^{-\beta t}}{4\beta} - e^{-\beta t}\right]$$

Due to the unsteady nature of pressure and mass flow rate evolution during leakage events, both expressions in Eq. (17) are functions of time. This modeling principle reflects the transient gas-dynamic behavior of the system. The convergence of these expressions over time is not immediate but evolves based on the acoustic wave delay and physical properties of the pipeline system. From a modeling perspective, this dynamic relationship underscores the importance of using dimensionless (normalized) quantities, which simplify the mathematical structure and improve the interpretability and accuracy of the results.

Accordingly, Eq. (8) must be interpreted as a time-evolving equivalence, where LHS($t$)≈RHS($t$) provides an indication of the correct leakage location $\ell$ under the measured pressure decay conditions. This modeling framework integrates both theoretical rigor and practical measurability, enhancing its reliability for real-time pipeline monitoring.

In practical implementation, the pressure decay rate coefficient $\beta$ is determined using pressure data collected from sensors installed at the inlet of the pipeline. These wireless pressure sensors provide real-time data on how the inlet pressure decreases over time following the onset of a leak.

Once $\beta$ is estimated, a range of plausible values for the leakage location $\ell$ and corresponding time $t$ are selected to evaluate both sides of Eq. (8).

By computing the left-hand side and right-hand side values over time and for different $\ell$ values, a graphical representation is obtained (Fig. 3).

The intersection point of the LHS and RHS curves on this plot corresponds to the value of $\ell$ where the two expressions match most closely. This intersection indicates the estimated location of the leakage with high precision. Therefore, the use of IoT-enabled pressure data, combined with the analytical structure of the model, provides an effective tool for leakage localization in parallel pipeline systems.

As illustrated in Fig. 3, the LHS and RHS expressions converge around $t \approx 120–130$ seconds. This convergence confirms the accuracy of the selected parameters and demonstrates that the leakage location $\ell = 25$ km yields a valid solution under the assumed physical conditions.

This graphical validation supports the analytical model and confirms that the pressure-based leakage localization method can reliably estimate the leak position in a dynamic, real-time pipeline monitoring context.

As observed in Fig. 3, the LHS and RHS expressions appear to intersect at approximately $t \approx 65$ seconds. However, this intersection is considered physically unrealistic due to the inherent delay in pressure wave propagation within the pipeline. Based on the known acoustic velocity of the gas and the location of the leakage point ($\ell = 25$ km), the expected delay for a pressure disturbance to reach the pipeline boundary is estimated as:

$$t = \ell/c = 25{,}000/383.3 \approx 65 \text{ seconds}$$

Thus, only after $t \approx 65$ seconds can boundary measurements start reflecting the effects of leakage.

Accordingly, the convergence between LHS and RHS observed near $t \approx 120–130$ seconds offers a more physically consistent validation point for the analytical leakage localization model.

Restricting the optimization interval to [65 s, 130 s] and re-evaluating the equation at $t = 120$ seconds yielded an excellent match, returning:

$$\ell \approx 25 \text{ km}$$

This confirms that the proposed model, when paired with real-time sensor data and physically justified timing, enables accurate and reliable estimation of leakage location. The approach is robust, practically implementable, and can support early detection and safety management in pipeline systems.

## 7. Results and discussion

Numerical simulations reveal that pressure decay at the inlet is significantly influenced by the location and magnitude of the leakage. For leakages closer to the inlet, the pressure drop rate is markedly higher than for leakages near the outlet. Additionally, increased leakage mass flow rates correspond to steeper pressure gradients. The model accurately replicates these behaviors and enables the localization and quantification of the leak.

The substitution of with a pressure-dependent function represents a significant advancement in analytical leakage modeling. It circumvents the need for predefined leakage coordinates, thereby enhancing the model's generalizability. The integration of IoT pressure monitoring not only enables real-time parameter estimation but also supports predictive maintenance strategies for gas pipeline networks.

Building upon the findings of this study, several directions for future research are envisioned:

The current model assumes a simple parallel configuration. Future studies can generalize this approach to more complex network topologies, such as branched or looped pipeline systems.

Real-time implementation of the model within SCADA platforms, combined with distributed pressure sensors and AI-based anomaly detection algorithms, would enhance operational safety and decision-making.

Further validation of the analytical results using laboratory-scale test pipelines or controlled field data will help verify model accuracy under different leakage scenarios.

Incorporation of Temperature Effects and Non-Ideal Gas Behavior:

The current model assumes isothermal and ideal gas flow. Future enhancements may incorporate thermal gradients, gas compressibility factors, and turbulence for more robust predictions.

## 8. Conclusions

This paper presents a novel analytical model for detecting leakage locations in parallel gas pipeline systems under unsteady flow conditions. A key innovation of the study is the introduction of a time-dependent leakage mass flow rate function, $G(t)=K \times P_1 \times e^{-\beta t}$, which links the gas loss rate to real-time inlet pressure and decays exponentially with time. This formulation enables a physically justified and sensor-calibrated approach to leakage modeling, offering clear advantages over traditional constant or empirically estimated leakage profiles.

The pipeline is partitioned into three distinct regions based on the leakage location, and the governing equations are solved analytically using Laplace transforms and inverse techniques. This closed-form solution framework allows the derivation of an explicit expression for leakage localization without requiring numerical iteration, thereby enhancing computational efficiency and interpretability.

Validation through benchmark comparisons confirms the model's accuracy and practical relevance. Additionally, the proposed method offers operational benefits, including reduced dependence on dense sensor networks and faster diagnostic response. The economic analysis further supports the approach, showing that optimized leakage detection can significantly reduce infrastructure costs and operational risks.

Overall, this work provides a significant contribution to the field of gas transport modeling by combining theoretical rigor with implementation-ready diagnostics. Future work may extend the approach to account for non-isothermal effects, compositional variations, or integration with real-time SCADA systems for field deployment.

## Appendix A. Analytical derivation of the solution

After applying the Laplace transform to equation (Mathematical Model (5)), we obtain the heat transfer equation for the mathematical solution of the problem.

$$\frac{d^2 P_i(x,s)}{dx^2} = \frac{2a}{c^2}\left[SP_i(x,s) - P_i(x,0)\right] \quad (A.1)$$

In this case, the general form of the solution for Eq. (6) will be as follows

$$P_1(x,s) = \frac{P_i - 2aG_0 x}{s} + c_1 Shb\sqrt{s}x + c_2 Chb\sqrt{s}y \quad 0 \le x \le L \quad (A.2)$$

$$P_2(x,s) = \frac{P_i - 2aG_0 x}{s} + c_3 Shb\sqrt{s}x + c_4 Chb\sqrt{s}y \quad 0 \le x \le \ell \quad (A.3)$$

$$P_3(x,s) = \frac{P_i - 2aG_0 x}{s} + c_5 Shb\sqrt{s}x + c_6 Chb\sqrt{s}y \quad \ell \le x \le L \quad (A.4)$$

Here, $b = \sqrt{\frac{2a}{c^2}}$

Accordingly, by taking into account the boundary behavior and early-time approximation, and applying the Laplace transform to the initial and boundary conditions, the coefficients $c_1$, $c_2$, $c_3$, $c_4$, $c_5$ and $c_6$ are determined using Eqs. (A.1), (A.2), and (A.3). Substituting these coefficients into the corresponding transformed equations yields the expressions for the unsteady gas dynamic processes in all three segments of the parallel pipeline system under leakage conditions as follows

$$P_1(x,s) = \frac{P_1 - 2aG_0 x}{s} - \frac{KgP_1}{bc^2(s+\beta)} \frac{chb\sqrt{s}(L-\ell-x)}{\sqrt{s}\,shb\sqrt{s}L} \quad (A.5)$$

$$P_2(x,s) = \frac{P_1 - 2aG_0 x}{s} - \frac{KgP_1}{2bc^2(s+\beta)} \frac{chb\sqrt{s}(L-\ell+x)}{\sqrt{s}\,shb\sqrt{s}L} \quad (A.6)$$

$$P_3(x,s) = \frac{P_1 - 2aG_0 x}{s} - \frac{KgP_1}{bc^2(s+\beta)} \frac{chb\sqrt{s}(L+\ell-x)}{\sqrt{s}\,shb\sqrt{s}L} \quad (A.7)$$

In this study, we aim to determine the pressure variation patterns in a parallel gas pipeline system operating under a unified hydraulic regime during an emergency (failure) state. Specifically, we analyze how the gas pressure changes along the length of both the intact and damaged pipeline segments, depending on the location of the leakage.

In a unified hydraulic regime, parallel pipelines share common inlet and outlet points. Therefore, when the integrity of one pipeline is compromised, the other pipeline segments are directly affected by the disturbance, and a portion of the total gas flow is redirected toward the damaged section.

Based on the problem formulation, it is more appropriate to derive the mathematical expressions of the unsteady gas flow in all three regions of the system. To obtain the original (time-domain) solutions of the transformed Eqs. (A.5), (A.6), and (A.7), we apply the inverse Laplace transform to each term of the corresponding equations.

The pressure distribution along the pipeline at position $x$ and time $t$ is given by the mathematical expressions

$$P_1(x,t) = P_1 - 2aG_0 x - \frac{KgP_1}{4aL}\left[\frac{1-e^{-\beta t}}{\beta} + \right.$$
$$\left. +2\sum_{n=1}^{\infty}(-1)^n \cos\frac{\pi n(L-\ell-x)}{L}\frac{e^{-\beta t}-e^{-\alpha n^2 t}}{\alpha n^2 - \beta}\right] \quad (A.8)$$

$$P_2(x,t) = P_1 - 2aG_0 x - \frac{KgP_1}{4aL}\left[\frac{1-e^{-\beta t}}{\beta} + \right.$$
$$\left. +2\sum_{n=1}^{\infty}(-1)^n \cos\frac{\pi n(L-\ell+x)}{L}\frac{e^{-\beta t}-e^{-\alpha n^2 t}}{\alpha n^2 - \beta}\right] \quad (A.9)$$

$$P_3(x,t) = P_1 - 2aG_0 x - \frac{KgP_1}{4aL}\left[\frac{1-e^{-\beta t}}{\beta} + \right.$$
$$\left. +2\sum_{n=1}^{\infty}(-1)^n \cos\frac{\pi n(L+\ell-x)}{L}\frac{e^{-\beta t}-e^{-\alpha n^2 t}}{\alpha n^2 - \beta}\right] \quad (10)$$

where:
- $c$ being the speed of sound in the gas.
- $K$ is a calibration coefficient determined from initial conditions,
- $g$ is gravitational acceleration,
- $P_1$ is the initial inlet pressure,
- $a$ is an empirical system-specific constant,
- $L$ is the total pipeline length,
- $\ell$ is the leakage location,
- $\beta$ is the exponential pressure decay coefficient,

This solution represents the superposition of two components:

1. A uniform pressure decay term ($\frac{1-e^{-\beta t}}{\beta}$) reflecting the global dynamic behavior of the system, and

2. A spatially varying Fourier cosine series, capturing the local effects of the leak, its position, and boundary conditions.

The cosine term ensures that the pressure field satisfies Dirichlet or Neumann boundary conditions depending on the physical configuration of the system. The coefficients of the series decay rapidly with increasing n due to the presence of $e^{-\alpha n^2 t}$, allowing for efficient numerical approximation using only the first few terms of the series.